# Catalan tableaux and the asymmetric exclusion process

Xavier Gérard Viennot
LaBRI, Université Bordeaux 1, 33405 Talence, France
viennot (à) labri (point) fr

**Abstract**
   The purpose of this paper is twofold. First we answer to a question asked by Steingrimsson and Williams about certain permutation tableaux: we construct a bijection between binary trees and the so-called Catalan tableaux. These tableaux are certain Ferrers (or Young) diagrams filled with some 0's and 1's, satisfying a certain hook condition, and are enumerated by the Catalan numbers. They form a subclass of the permutation tableaux, enumerated by n!, introduced by Postnikov in his study of totally non negative Grassmannians and networks.
   Secondly we relate this new Catalan bijection with the totally asymmetric exclusion process(TASEP), a very rich and well studied 1D gas model in statistical mechanics of non-equilibrium systems. We continue some combinatorial understanding of that model, in the spirit of works by Shapiro, Zeilberger and more recently by Brak, Essam, Rechnitzer, Corteel, Williams, Duchi and Schaeffer. Emphasis is made on the non-classical notion of canopy of a binary tree, analog of the classical up-down sequence of a permutation.

**Résumé**
   Le but de cet article est double. Tout d'abord nous répondons à une question posée par Steingrimsson et Williams sur certains "tableaux de permutation": nous construisons une bijection entre les arbres binaires et la classe d'objets appelés tableaux de Catalan. Ces tableaux sont des diagrammes de Ferrers remplis avec des 0 et des 1, soumis à une certaine condition d'équerre et sont énumérés par les nombres de Catalan. Ils forment une sous-classe des tableaux de permutation, énumérés par n!, objets introduits par Postnikov dans son étude des Grassmanniennes totalement positives et des "réseaux".
   Deuxièmement, nous montrons le lien entre cette nouvelle bijection du monde Catalan avec le modèle de gaz d'exclusion asymétrique (TASEP), un modèle sur un réseau à une dimension, très riche et très étudié en physique statistique des systèmes hors d'équilibre. Nous continuons l'étude pour une compréhension combinatoire du modèle, dans la lignée des travaux de Shapiro et Zeilberger, et plus récemment de Brak, Essam, Rechnitzer, Corteel,Williams, Duchi et Schaeffer. Une notion centrale est la notion peu classique de canopée d'un arbre binaire, analogue de celle bien classique de forme d'une permutation (succession des montées et des descentes).

## §1 Introduction

   Following work of Postnikov about totally non negative Grassmannians [14], Einar Steingrimsson and Lauren Williams introduced in [16] the notion of permutation tableaux as Ferrers diagrams filled with 0 and 1 and satisfying a certain hook condition. These tableaux are in bijection with permutations. If each column has only one 1, then such tableaux are enumerated by the Catalan numbers. We propose to call them *Catalan tableaux*. In their seminal paper, Steingrimsson and Williams exhibit a bijection between these tableaux and a class of permutations defined by a certain forbidden pattern. They ask the question of finding a bijection between these tableaux and one of the classical objects enumerated by Catalan numbers. In this paper we construct a bijection between Catalan tableaux and binary trees. Other solutions with other Catalan objects have been given by Burstein [3] with non-crossing partitions and also by Corteel, Eriksen and Riefegerste (quoted in [16]).

   In the second part of this paper, we use this new Catalan bijection for unifying and going further in the combinatorial study of the well known model in physics called totally asymmetric exclusion process (TASEP). This model is the archetype of non-equilibrium gas models, with explicit resolution, deep properties and existence of phase transitions. In this model, particles are hopping to the right in a one-dimensional strip of $n$ cells. At any time each particle can jump to next cell with probability $dt$ (if it is not occupied), enter the strip on the left with probability $\alpha dt$ and leave the strip on the right with probability $\beta dt$. The exclusion condition means that there is at most one particle in each of the $n$ cells. The model can also be described by a discrete Markov chain on the set of the $2^n$ possible states formed by words $w$ of length $n$ on an alphabet with 2 letters ("occupied" or "empty" cell). There exist unique stationary probabilities $prob(w;\alpha,\beta)$ for any state $w$.

Many works have been done by physicists, in particular Derrida and coautors. A method with a *matrix ansatz* has been introduced by Derrida et al. [9] and appears to be a very powerful tool for giving various different expressions for the probabilites $prob(w;\alpha,\beta)$. Surprisingly, Catalan numbers appear and a very rich and deep combinatorial understanding has been developped, in particular by Shapiro, Zeilberger [15], Brak, Essam [2], Corteel, Brak, Rechnitzer, Essam [4], Corteel, Williams [5], involving Catalan combinatorics of weigthed paths and Catalan tableaux, in relation with different solutions of the matrix ansatz. In fact the papers [2], [4], [5] deal with the more general model ASEP (or PASEP for partially asymmetric exclusion process), when the particles are allowed to move backward to the left with probabilities $qdt$.

In another spirit, Duchi and Schaeffer [10], avoiding the matrix ansatz, give a direct combinatorial interpretation of the appearance of the Catalan numbers by constructing a larger Markov chain with uniform stationary distribution which projects on the Markov chain of the TASEP. An analog philosophy is developed in the paper [4] for a more general model.

In this paper we use our new Catalan bijection between tableaux and binary trees to go further in the combinatorial understanding of the TASEP. Combining that bijection with some other classical bijections from the "Catalan garden", we relate two different interpretations of the stationary distribution of the TASEP model: the one given by Shapiro, Zeilberger [15] (for the parameters $\alpha = \beta = 1$) in term of pair of paths (also called *staircase polygons* or *parallelogram polyominoes*) and the interpretation given by Corteel, Williams [5] in term of Catalan tableaux. This leads us to give a direct bijection, without going through binary trees, between Catalan tableaux and pair of paths ($\omega, \eta$) sending the shape of the tableau (corresponding to the state $w$ of the TASEP) to the path $\omega$.

The combinatorics we develop here is strongly related to the notion of *canopy* of a binary tree introduced by the author for the Lascouxfest [17]. Our main bijection associates the shape of the Catalan tableau to the canopy of the binary tree. Combining that bijection with the combinatorics developed in [17] reveal deeper combinatorics hidden behind this simple TASEP model. We give an interpretation of the stationary probabilities with general parameter $\alpha$ and $\beta$ in term of binary trees, showing in an obvious way the symmetry of the TASEP when exchanging the parameters $\alpha$ and $\beta$. Also we relate the Catalan tableaux interpretation of Corteel, Williams [5] with that given by Duchi, Schaeffer [10] in term of pair of paths (with parameter $\alpha$ and $\beta$).

Another interest of the bijection between binary trees and Catalan tableaux relies in the relation with the the Loday-Ronco Hopf algebra of binary trees [12]. We can define directly in term of Catalan tableaux the product and the coproduct of the dual Hopf algebra (paper in preparation with Aval [1]).

**§2 Catalan Tableaux**

Recall that a partition $\lambda = (\lambda_1, ..., \lambda_k)$ is a weakly decreasing sequence of non-negative intergers. We denote by $F_\lambda$ the *Ferrers diagram* (also called *Young diagram*) associated to the partition $\lambda$. It is a diagram having $m$ cells, with $\lambda_i$ cells in the i-th row. For convenience with the bijections involved in this paper, we will use neither of the classical ways of drawing Ferrers diagrams ("french" or "english" way). The diagram will be right-justified, and the rows written from bottom to top as in the "french" way, see Figure 1.

**Definition 2.1.** A *Catalan tableau* is a partition $\lambda$ such that the diagram $F_\lambda$ is contained in a $k \times (n-k)$ rectangle ($k$ rows and ($n$-$k$) columns), together with a filling of the cells with 0's and 1's satisfiying the two following properties:
(*i*) each column of the rectangle contains exactly one 1,
(*ii*) there is no 0 which has a 1 below it in the same column and to its right in the same row.
The integer $n$ is called the *index* of the tableau.

Catalan tableaux are displayed on Figures 1 and 4 (top). From condition (*i*) the Ferrers diagram must have (*n-k*) columns, whereas the number of rows may be smaller than $k$. For a fixed integer *n*, we know from Steigrimsson and Williams [16], that the number of such Catalan tableaux is the Catalan number $C_n$. We describe now a bijection between Catalan tableaux and binary trees, thus solving a question asked by the two authors in their paper.

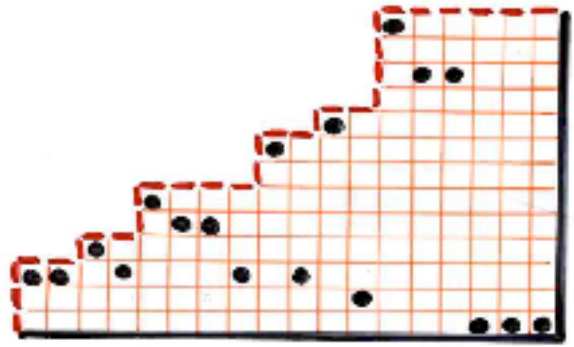

FIGURE 1. A Catalan tableau
(the 0's and 1's have been replaced by empty cell or black dot)

Recall that a *binary tree* $B$ is classically defined recursively as being: else the empty binary tree, else a triple (*L,r,R*) where $r$ is the *root* of the binary tree, and *L* (resp. *R*) is a binary tree called the *right* (resp. *left*) *subtree* rooted at the vertex $r$. The number of binary trees having $n$ vertices is the Catalan number $C_n$. We will draw binary trees with the root at the bottom, North-West (NW) (resp. North-East (NE)) edges will correspond to the link relating a vertex to its *left* (resp. *right*) *son*, i.e. the root of its left (resp. right) subtree. The *right* (resp. *left*) *branch* of $B$ is the sequence of vertices starting from the root and going from one vertex to another by taking the right (resp. left) son of each *father*.

In order to describe our bijection $\phi$ from Catalan tableaux to binary trees, we need to introduce a combinatorial object interpolating between binary trees and paths. We define a *tailed binary tree* as to be a pair $Q = (B, \omega)$ formed by a binary tree $B$ together with a path $\omega$ having NE and SE steps, starting at one of the vertices of the right branch of the binary tree (for $B$ non empty). Moreover, if the path is not empty, the first elementary step is SE. The path $\omega$ is called the *tail* of $Q$. Remark that usually an empty path is reduced to a point; in that case we should take the convention that this unique vertex is the last vertex of the right branch of $B$. Such objects are visualized on Figure 2.

The *profile* $\pi(T)$ of a Catalan tableau $T$ is the sequence of North (N) and East (E) steps of the NW border of the tableau, starting from the SW corner of the tableau, with the first step (always N) deleted.

**From Catalan tableaux to binary trees: the algorithm** $T \to B = \phi(T)$

We start from a Catalan tableau $T$ (displayed on Figure 2), rotated clockwise 45° and we draw its NW border as a path $\omega$ with NE (in red) and SE (in blue) steps. The SE (resp. NE) edges of the path are going to become the left and right edges of the binary tree $\phi(T)$. We define a *corner* of the path $\omega$ as to be a pair of successive edges NE and SE.

**algorithm 2.2,** described through steps (a),...,(e).

(a) If the path $\omega$ has no corner, i.e. the path is reduced to a sequence of $n$ NE steps, then the binary tree $B = \phi(T)$ is reduced to a sequence of $n$ vertices connected by only right edges, in other words $B$ is the path $\omega$ with its first step deleted. The algorithm stops.

(b) Else, we take the first corner of the profile. The path is splitted into $\omega = (s_0,...,s_k,s_{k+1},w)$ where the edges $(s_0,s_1),....,(s_{k-1},s_k)$ are NE (in red) while the edge $(s_k,s_{k+1})$ is SE (in blue). We have denoted by $w$ the remaining part of the path $\omega$ after the vertex $s_k$. The algorithm starts with $Q = (B,v)$ with $B$ being the binary tree reduced to its right branch which is the sequence of vertices $(s_0,...,s_k)$. The tail $v$ of $Q = (B,v)$ is the path $(s_k,s_{k+1},w)$ starting at the vertex $s_k$.

In general we have a pair $Q = (B,v)$ where $B$ is a binary tree having $(s_0,...,s_k)$ as right branch and $v$ is a path starting with a SE step (blue) at vertex $s_r$ with $0 \leq r \leq k$.

(c) If there is a 0 in the first corner of the profile $\pi(T)$, then the tableau $T$ is replaced by the tableau obtained from $T$ by deleting the row starting form that corner. The pair $Q = (B,v)$ is replaced by $Q' = (B,v')$ where $v'$ is the same path as $v$ considered as a sequence of SE or NE edges, but now starting at the vertex $s_{r-1}$ of the binary tree $B$, father of the vertex $s_r$.

(d) Else, there is a 1 in the first corner of the profile $\pi(T)$. The tableau $T$ is replaced by the tableau obtained from $T$ by deleting the column starting from that corner. The pair $Q = (B,v)$ is replaced by the pair $Q' = (B',v')$ defined as follows.

We split the path $v$ into three parts: its first step (SE) followed by a (possibly empty) sequence of NE steps, followed by a path $v'$ possibly empty, or starting by a SE step. We write $v = (t_0,t_1,....,t_p,...,t_q)$ where the fisrt step $(t_0,t_1)$ is SE, the sequence $\tau = (t_1,...,t_p)$ is a path having only NE steps, and $v' = (t_p,...,t_q)$ starts with a SE step or is

*4*

empty $(p=q)$. In our notation $t_0 = s_r$.

Then the binary tree $B'$ is obtained from $B$ by replacing the subtree $B_{s_r}$ rooted at the vertex $s_r$ by the binary tree $(B_{s_r}, s_r, \tau)$. We have considered the sequence $\tau$ as a binary tree reduced to its right branch. The path $v'$ is now starting at the last vertex of the right branch of $B'$, that is the last vertex $t_p$ of the sequence $\tau$.

(e)  Repeat step (c) untill one get $Q = (B,v)$ with $v = \emptyset$. Then the binary tree $B$ has always an empty left subtree. The binary tree $\phi(T)$ is obtained from $B$ by deleting the right edge starting from the root.
**end of algorithm 2.2**

An example of this recursive algorithm is displayed on Figure 2. At each step we have colored in yellow the row or the column which is going to be deleted from the Catalan tableau $T$. At each deletion, the remaining tableau is again a Catalan tableau. At step (c), the corner choosen has always a cell filled with a 1 (or black dot) below in the same column. At each step of the algorithm, the profile of the running tableau is (up to a 45° rotation) the sequence formed by the right branch of $B$ followed by the tail $v$ of $Q = (B,v)$.

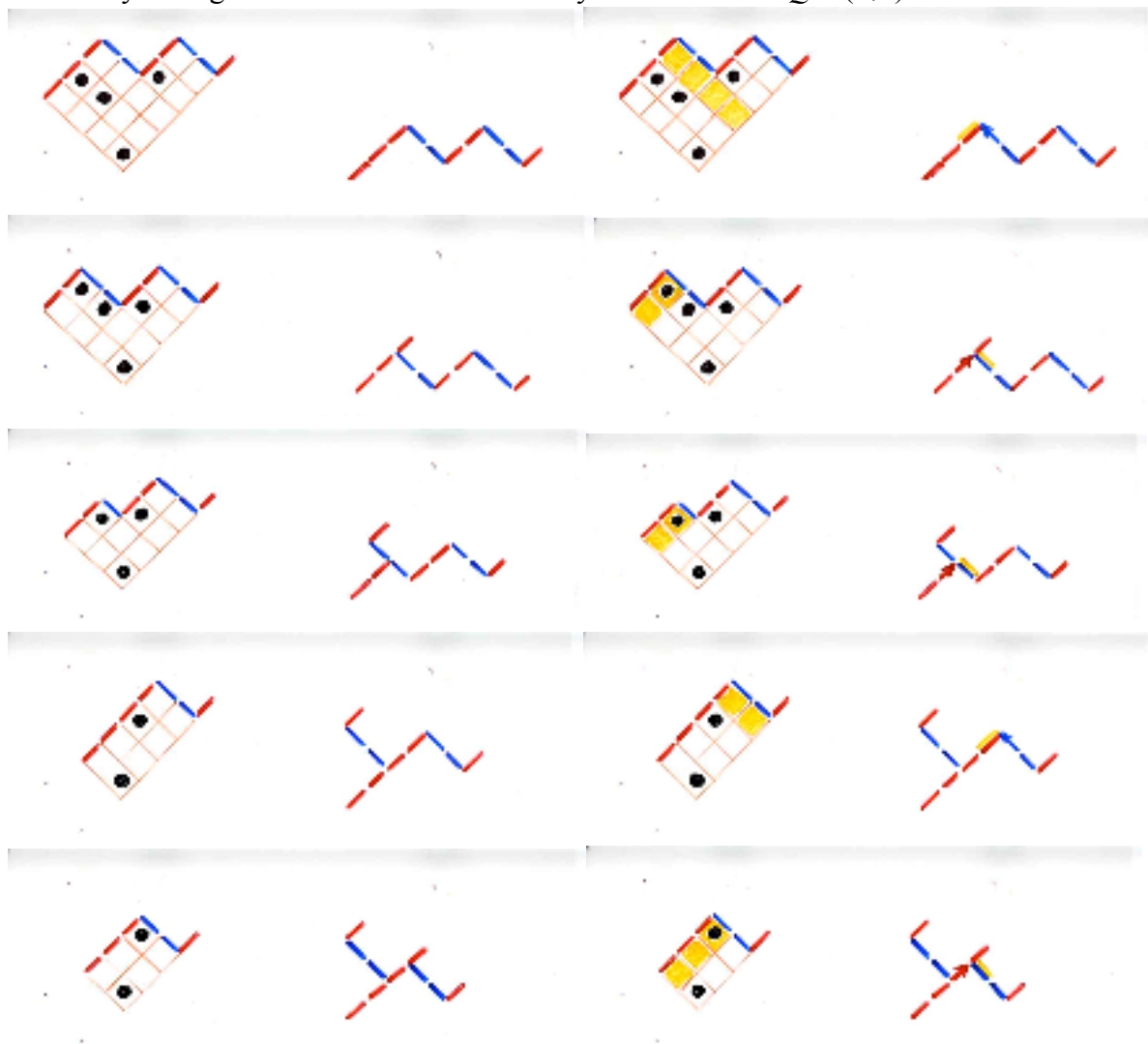

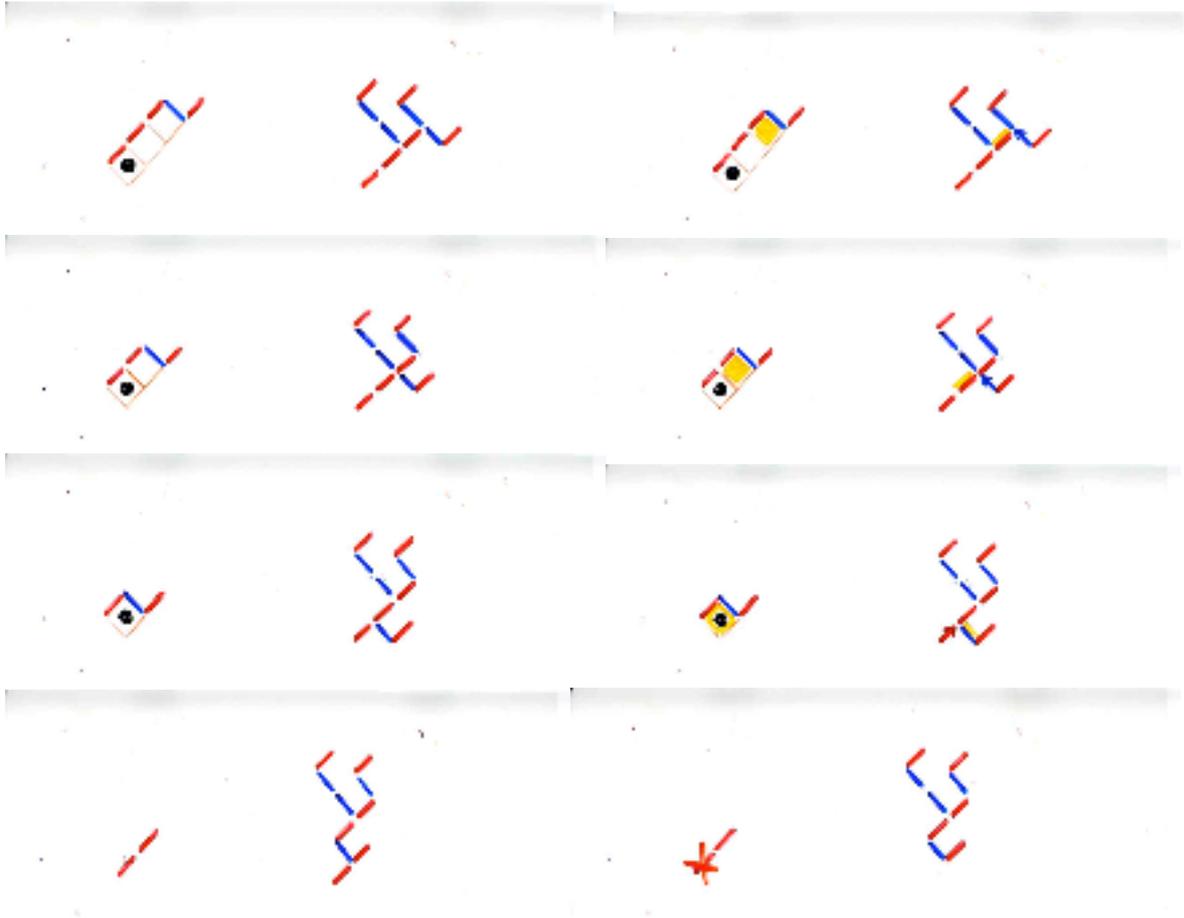

FIGURE 2. algorithm 2.2: the bijection $\phi$ between Catalan tableaux and binary trees

**Proposition 2.3.** *The map $\phi:T \to B$ defined by algorithm 2.2 described above with steps (a),(b),(c),(d),(e) is a bijection between the set of Catalan tableaux having index $n$ and the set of binary trees with $n$ vertices.*

The proof can be done recursively by constructing the reverse bijection. This reverse bijection can also be deduced from the construction in §4 below. A natural question is to ask to which parameter in the binary tree corresponds the number $(n-k)$ of columns of the Catalan tableau $T$, or to which parameter in the tableau corresponds the *length* (= number of vetices) of the right and the left branch in the binary tree $B$. What is the interpretation in the binary tree $B$ of the shape of the tableau $T$? The answers to these questions are given below. We introduce the following definitions.

Let $B$ be a binary tree with $n$ vertices. A classical notion is the notion of *symmetric order* (also called *inorder*) which is a total order on the set of vertices of a binary tree. It is defined recursively by the following: first follow the left subtree, then the root, then the right subtree. A classical bijection $\gamma$ is from binary trees with $n$ vertices to complete binary trees with $2n+1$ vertices, i.e. binary trees with no vertex having only one son. The tree $\gamma(B)$ is obtained form $B$ by adding an external vertex on the left (resp. right) of every vertex not having a left (resp. right) son (see Figure 4, top).

**Definition 2.4.** The *canopy* of a binary tree $B$ with $n$ vertices is the sequence $u = u_1...u_{n-1}$ of length ($n$-$1$) written with two letters $a$ (red) and $b$ (blue) such that the ith letter is $a$ (resp. $b$) iff the *ith* vertex of $B$ for the symmetric order has (resp. does not have) a right son.

The folowing definition, introduced by Dulucq and Guibert [11] is equivalent: the canopy $u$ of the binary tree $B$ is the word obtained by reading the leaves of the complete binary tree $\gamma(B)$ in symmetric order and writing the letter $a$ (resp.$b$) iff the leaf is a left (resp. right) daughter (=son) of her father. The fisrt and last leaf are not considered (giving always respectively an "$a$" and a "$b$").

**Proposition 2.5** *Let T be a Catalan tableau of index $n$ with profile $u = \pi(T)$ and let $v$ be the canopy of the binary tree $B = \phi(T)$ associated by the bijection defined by algorithm 2.2. Then the two words $u$ and $v$ are equal by identifying letters $a$ with N, and $b$ with S.*

The proof relies on an extension of the notion of symmetric order of a binary tree to tailed binary trees and by showing that this symmetric order is invariant through algorithm 2.2. Following Corteel and Williams [5], we say that an entry in a column is *restricted* if that entry is a 0 which lies above some 1. We say that a row is *unrestricted* if it does not contain a restricted entry. It is not difficult to prove the following characterisation of the length of the branches of the binary tree $B = \phi(T)$ :

**Proposition 2.6**. *Let T be a Catalan tableau and $B = \phi(T)$ be the associated binary tree. The length of the left branch (resp. right branch) is the number of 1's in the first row (resp. number of unrestricted rows) of the tableau T, augmented by 1.*

§3 **Relation with the TASEP**

The *totally asymmetric exclusion process* is a classical and intensively studied gas model in statistical mechanics. It is the archetype of non-equilibrium models, being simple enough to lead to explicit resolutions, but being surprisingly rich, with deep propreties and having phase transition of first and second order. It has appeared as a modelisation in various domains (biology, traffic flow, statistical physics, ...) and is connected to the theory of orthogonal polynomials. Recently some deep and fruitful connections has been found with combinatorics [2],[4],[5],[10],[15].

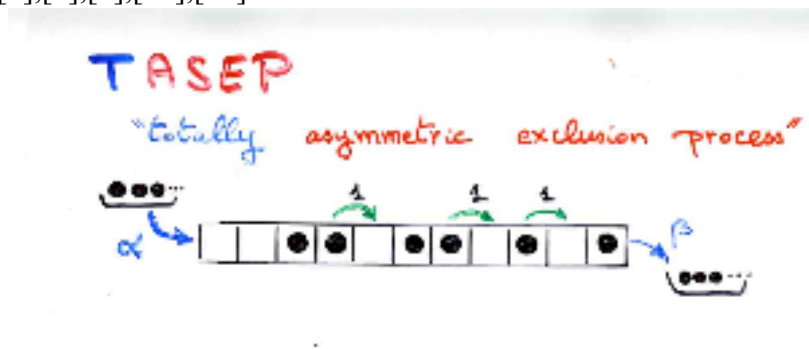

FIGURE 3. Illustration of the TASEP

The TASEP model describes a system of particles hopping to the right on a one-dimensional lattice of $n$ sites with the exclusion property, i.e. at most one particle is found on a site (or cell). In the continuous model, roughly speaking, particles may hop at any time to the right with rate $dt$ (if the cell is not occupied), may enter at the left of the strip with rate $\alpha dt$ and may exit the system at the right of the strip with rate $\beta dt$ (if possible).

The discrete model can be described as a Markov chain on the set of TASEP configurations, i.e. the $2^n$ words of length $n$ formed by the sequence of states ("occupied" or "empty") of the $n$ cells of the strip. At discrete time $t=0,1,..,i,..$ a wall (separation between two cells, including the borders of the strip) is choosen at random with uniform distribution. Transitions are made (if possible) by a particle traversing the wall with probability 1, $\alpha$, $\beta$, according to the three possible cases: respectively internal, leftmost, rightmost wall.

Both models leads to the same stationary distribution in the limit ($t \to \infty$) regime (see for example [7], [8] or [10] for a more detailed description). We encode a state $u$ by a path $\omega$ on the square lattice. It is the sequence of elementary steps North (N) and East (E) obtained by replacing in the sequence $u$ "cell occupied" (resp."empty") by the step N (resp.E). In the case $\alpha = \beta = 1$, we denote by $prob(u)$ the stationary probability associated to the state $u$.

One of the first explicit expression for that stationary distribution was given in 1982 by Shapiro and Zeilberger [15], with motivation from molecular biology.

(3.1) $$prob(u) = \frac{1}{C_{n+1}} \textit{ (number of paths } \eta \textit{ )},$$

where $\eta$ is a path on the square lattice, formed by elementary steps North or South, having the same endoints as the path $\omega$ associated to $u$, and located below the path $\omega$

. Note that "below" means that the path $\eta$ can touch the path $\omega$ but never cross it (see Figure 4). For a given path $\omega$ the number of such path $\eta$ has been given in 1955 by Narayana [13] as a determinant (extended by Kreweras). If $\lambda = (\lambda_1,...,\lambda_k)$ is the partition associated to the Ferrers diagram defined by the path $\omega$, then the number of paths $\eta$ in equation (3.1) is the determinant $\det\binom{\lambda_i + 1}{j - i + 1}_{1 \leq i, j \leq k}$.

A second expression in terms of Catalan tableaux has been given by Corteel and Williams [5]:

(3.2) $$prob(u) = \frac{1}{C_{n+1}} \textit{ (number of Catalan tableaux T)},$$

where $T$ is a Catalan tableau having as profile $\pi(T)$ the path $\omega$ associated to $u$.

Both expressions (3.1) and (3.2) can be deduced from the *ansatz matrix* method introduced by Derrida et al in [9]. Each identity corresponds to different pairs (D,E) of matrices satisfying the commutation rule $DE = D + E$. A direct combinatorial proof of (3.1) has been given by Duchi and Schaeffer [10] by constructing a Markov chain on the set of pair of paths ($\omega$, $\eta$), with uniform stationary distribution, and "projecting" on the TASEP

Markov chain. In the same spirit, direct combinatorial proof of (3.2) has been given by Corteel and Williams in [5] (in fact as a consequence of a more general construction involving the so-called permutation tableaux, enumerated by $n!$, and related to the PASEP model that is the model as above but where the particles can hop on the left with probability $qdt$.

A natural question is to relate combinatorially expression (3.1) and (3.2), that is to give a bijection from pair of paths $(\omega, \eta)$ to Catalan tableaux $T$ such that the profile of $T$ is the path $\omega$. Such correspondence, using the bijection $\phi$ constructed in section 2 is given in the next section.

## §4 Combinatorial relation between the interpretations of the stationary probabilities of the TASEP in term of paths and Catalan tableaux

This relation is visualized on Figure 4, combining the bijection between Catalan tableaux and binary trees with some classical bijections of the "Catalan garden". From a Catalan tableau $T$ of index $n+1$ (on Figure 4, $n=9$), we associate the binary tree $B = \phi(T)$. A bijection associating the pair of paths $(\omega, \eta)$ to the binary tree $B$ has been described and studied in Viennot [17] (see also Delest, Viennot [6]). It is summarized on Figure 4 and can be obtained by combining classical bijections of the "Catalan garden". We recall these constructions.

First, we associate a Dyck path $\xi$ to the binary tree $B$, by extending the binary tree $B$ to the complete binary tree $\gamma(B)$ and we read that tree in reverse prefix order (follow recursively fisrt the root, then the right subtree, then the left subtree). For each internal (resp. external = leaf) vertex, we put a NW (resp. SW elementary step), reading the Dyck path $\xi$ in a reverse way from right to left. Then a *parallelogram polyomino* $P$ (= *staircase polygon*) is associated to the Dyck path $\xi$. We read, as usual, the Dyck path $\xi$ from left to right. A *peak* (resp. *valley*) of the Dyck path is a succesion of NE-SE (resp. SE-NE) steps. To each peak having *height* $h$, we associate a vertical strip of $h$ cells. These strips are glued together in order to form the columns of the polyomino $P$ in a unique way so that the number of common vertical edges between two consecutive columns is the height of the valley located between the two corresponding peaks, augmented by 1. Now, by removing the two "corners" of the polyomino $P$ and by sliding by an elementary NW step the inferior border of $P$, we get the pair of paths $(\omega, \eta)$. It is a remarkable fact that the path $\omega$ is the canopy (up to a change of letters $a,b$ to N,S) of the binary tree $B$ (which is also the profile $\pi(T)$).

A second equivalent description of the bijection $B \rightarrow (\omega, \eta)$ is the following. Each left edge in $B$ (colored in blue on Figure 4) connects two vertices: a father and its son having the same "*right height*", that is the number of "*right edges*" (in red on Figure 4) on the unique path going from the root of $B$ to one of these two vertices. We totally order these left edges by the induced order of the associated fathers for the symmetric order of the binary tree $B$. (labeled $a,b,c,d$ on Figure 4).

We define the path $\omega$ as to be the path corresponding to the canopy of $B$ (with the correspondance $N = a, S = b$). The number of East step (blue) is the same as the number of left edges in the binary tree $B$. We number these steps from left to right. For each East step of $\omega$, we put below in the same column an horizontal step (colored purple on Figure 4) at a distance which is precisely the right height of the corresponding left edge of the binary tree $B$ (ordered $a,b,c,..$ as above). It is a combinatorial property of binary trees that the successive levels of such edges are non-decreasing, and thus can be joined in a unique way by vertical steps forming the path $\eta$ located below the path $\omega$ and having the same ending points. This construction is equivalent to the first one going through Dyck paths and polyominoes.

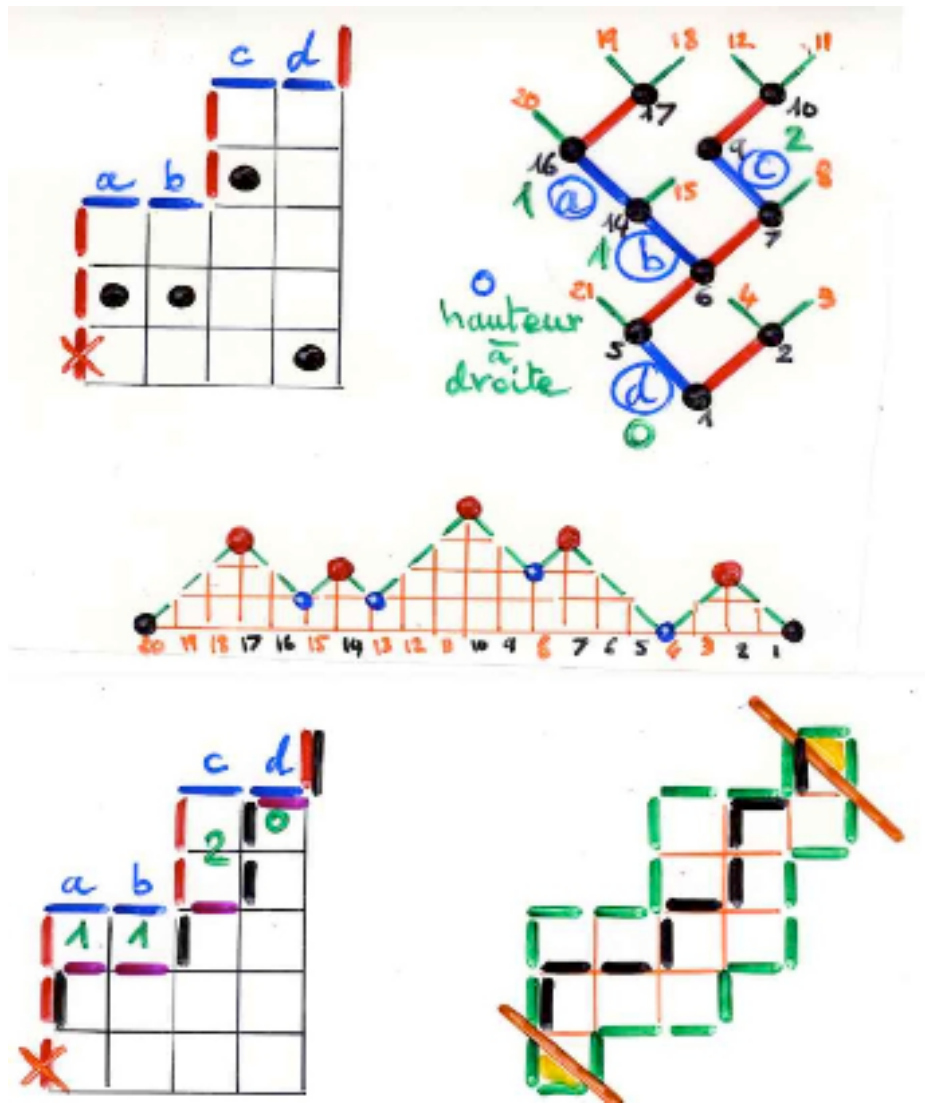

FIGURE 4. Catalan bijections $T \to B = \phi(T) \to \xi \to P \to (\omega, \eta)$.

We deduce the equivalence of the two formulae (3.1) and (3.2) for the stationary probabilities of the TASEP. As a consequence of the previous combinatorial constructions, we get a bijection from Catalan tableaux $T$ to pair of paths $(\omega, \eta)$. The description with right height of left edges in the binary tree enables us to give a direct definition of the bijection $T \to (\omega, \eta)$.

Let $T$ be a Catalan tableau. For each cell with a 1, in a first step we select all the cells above it (that is the *restricted* cells), and in a second step we select all the cells at the right of the selected cells in the first step. All these cells are colored yellow on Figure 5. For each column of the tableau $T$, we draw an horizontal edge (colored in blue on Figure 5) at a distance from the upper border which is equal to the number of white cells in the column, i.e. empty

*10*

cells not colored in yellow. These selected edges colored in blue define the second path $\eta$.

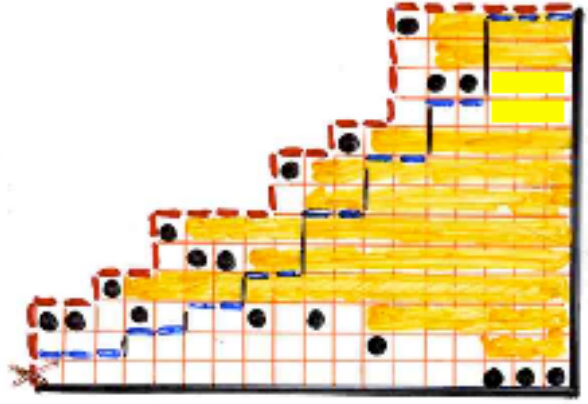

FIGURE 5. The direct bijection $T \to (\omega, \eta)$

## §5 The general TASEP with probabilities $\alpha$ and $\beta$.

We sumarized further combinatorial developments resulting from the main bijection $\phi$. The matrix ansatz of Derrida and al. [9] gives different interpretations of the stationary probabilities with parameters $\alpha$ and $\beta$ for the transitions probabilities entering or exiting the strip of cells of the TASEP. We denote by $prob(u; \alpha, \beta)$ the stationary probability for the state $u$ of length $n$ with parameters $\alpha$ and $\beta$.

In [5] Corteel and Williams refer to the matrix ansatz [9] with a certain pair of matrices $(D, E)$ and vectors $(V, W)$ to deduce an interpretation of $prob(u; \alpha, \beta)$ in terms of Catalan tableaux. with the parameters "number of *unrestricted* rows" and "number of 1's in the first row" of $T$ (see theorem 3.1 of [5] in the case $q = 0$). From proposition 2.6, we deduce the following interpretation:

(5.1) $$prob(u; \alpha, \beta) = \frac{1}{Z_n} \sum_B p(B),$$

where the summation is over all binary trees $B$ having canopy $u$, the weight $p(B)$ is defined by $p(B) = \alpha^{-lb(B)} \beta^{-rb(B)}$ where $lb(B)$ (resp. $rb(B)$) denote the length of the left (resp. the right) branch diminished by 1. The normalisation factor $Z_n$ is defined as usual by:

(5.2) $$Z_n = \sum_B p(B),$$

where the summation is over all binary trees with $(n+1)$ vertices.

Formula (5.1) is symmetric when exchanging $\alpha$ and $\beta$, left and right. We recover in an obvious way the duality of the TASEP exchanging left and right, empty and occupied sites. Using the bijections described in section 4, we can "transport" this interpretation in terms of pair of paths $(\omega, \eta)$.

(5.3) $$prob(u; \alpha, \beta) = \frac{1}{Z_n} \sum_{(\omega, \eta)} p(\omega, \eta),$$

where the summation is over all pair of paths $(\omega, \eta)$ of length $n$ and $p(\omega, \eta)$ is defined by $p(\omega, \eta) = \alpha^{-f(\omega,\eta)} \beta^{-g(\omega,\eta)}$, where $f(\omega, \eta)$ is the number of *contacts* between the two paths along horizontal edges and $g(\omega, \eta)$ is the number of steps N at the end of the path $\eta$.

Such pair $(\omega, \eta)$ can be put in bijection with the so-called *2-colored Motzkin paths*, or can be viewed as the configurations introduced by Duchi and Schaeffer in [10]. The interpretation given here is different from the interpratation given in their paper with theorem 3.3. But we can relate both interpretations by using a non classical bijection of the "Catalan garden" among binary trees transforming the symmetric order into the prefix order.

Moreover, from Proposition 2.6 and using some combinatorics of binary trees and *ballot numbers*, it is easy to prove combinatorially the known formula giving $Z_n$ (see for example Derrida [7]):

$$(5.4) \qquad Z_n = \sum_{i=1}^{n} \frac{i}{2n-i} \binom{2n-i}{n} \frac{\alpha^{-(i+1)} - \beta^{-(i+1)}}{\alpha^{-1} - \beta^{-1}}.$$

**Final remark.**

The combinatorics of the bijection $\phi$ is related to the study in Viennot [17] with the "*jeu de taquin*" for binary trees, and is related to the Loday-Ronco Hopf algebra of binary trees [12]. In Aval, Viennot [1] we describe the product and the co-product of the dual Hopf algebra directly in term of Catalan tableaux.. The computations become particularly simple.

**Acknowledgements**

Many thanks to Lauren Williams for fruitful discussions about the combinatorics of Catalan tableaux and related topics.